\documentclass{jnmp01b}

\usepackage{amsmath}

\renewcommand{\thepage}{R\arabic{page}}

\numberwithin{equation}{section}

\DeclareMathOperator{\Log}{Log}

\allowdisplaybreaks

\makeatletter
\def\cases#1{\left\{\,\vcenter{\normalbaselines\m@th
    \ialign{$##\hfil$&\quad##\hfil\crcr#1\crcr}}\right.}
\def\matrix#1{\null\,\vcenter{\normalbaselines\m@th
    \ialign{\hfil$##$\hfil&&\quad\hfil$##$\hfil\crcr
      \mathstrut\crcr\noalign{\kern-\baselineskip}
      #1\crcr\mathstrut\crcr\noalign{\kern-\baselineskip}}}\,}
\makeatother

\def\dps{\displaystyle}
\def\N{{\bf {N}}}
\def\Z{{\bf {Z}}}

\newenvironment{entry}[1]%
{\begin{trivlist}\item[]\textbf{Entry #1.}\ \itshape}
{\end{trivlist}}
\newenvironment{proposition}[1]%
{\begin{trivlist}\item[]\textbf{Proposition #1.}\ \itshape}
{\end{trivlist}}
\newenvironment{corollary}[1]%
{\begin{trivlist}\item[]\textbf{Corollary #1.}\ \itshape}
{\end{trivlist}}
\newenvironment{theorem}[1]%
{\begin{trivlist}\item[]\textbf{Theorem #1.}\ \itshape}
{\end{trivlist}}

\newenvironment{remark}[1]%
{\begin{trivlist}\item[]\textbf{Remark #1.}\ }
{\end{trivlist}}

\makeatletter
\DeclareRobustCommand{\primfrac}[1]{%
  \PackageWarning{amsmath}{%
Foreign command \@backslashchar#1; %
\protect\frac\space or \protect\genfrac\space should be used instead%
  }
  \global\@xp\let\csname#1\@xp\endcsname\csname @@#1\endcsname
  \csname#1\endcsname
}
\makeatother

\begin{document}

\renewcommand{\evenhead}{... by B A\ Kupershmidt}
\renewcommand{\oddhead}{Book Review}


\thispagestyle{empty}

\begin{flushleft}
\footnotesize \sf
Journal of Nonlinear Mathematical Physics \qquad 2000, V.7, N~2,
\pageref{firstpage}--\pageref{lastpage}.
\hfill {\sc Book Review}
\end{flushleft}

\vspace{-5mm}


\Name{Book Review by B A\ Kupershmidt}

\label{firstpage}

\vspace{-2mm}

\Adress{The University of Tennessee Space Institute,
 USA}


\vskip 12pt

\noindent
Five books are reviewed, namely

\vskip 10pt

\noindent {\bf Bruce C\ Berndt}:
{\it Ramanujan's Notebooks.\ Part I}.\/ (With a
foreword by {\bf S Chandrasekhar}). Springer-Verlag, New York Berlin, 1985.
357 pages. 

\vskip 10pt

\noindent{\bf ---}:
{\it Ramanujan's Notebooks.\ Part II}.\/
Springer-Verlag, New York Berlin, 1989. 359
pages. 

\vskip 10pt

\noindent{\bf ---}:
{\it Ramanujan's Notebooks.\ Part III}.\/
Springer-Verlag, New York, 1991. 510 pages.

\vskip 10pt

\noindent{\bf ---}:
{\it Ramanujan's Notebooks.\ Part IV}.\/
Springer-Verlag, New York, 1994. 451 pages.

\vskip 10pt

\noindent{\bf ---}:
{\it Ramanujan's Notebooks.\ Part V}.\/
Springer-Verlag, New York, 1998. 624 pages.

\vspace{7mm}

\begin{flushleft}
\large \bf \ldots And Free Lunch for All.\\
A Review of Bruce C Berndt's \textit{Ramanujan's Notebooks}, Parts I -- V.
\end{flushleft}

\vspace{1mm}

\begin{flushright}
\parbox{6cm}{\sffamily \small What Mozart was to music and Einstein 
was to physics, Ramanujan was to math. \\
\null\qquad Clifford Stoll}
\end{flushright}

\medskip

\noindent
Let $p(n) $ denote the number of partitions of a positive integer $n$.  Anyone who can dream 
up the formulae
\begin{gather*}
p(4) + p(9) x + p(14) x^2 + p(19) x^3 + p (24) x^4 + ... 
\\
\qquad = 5 {\{(1 - x^5) (1 - x^{10}) (1 - x^{15}) (1 - x^{20}) ...\}^5 \over 
\{(1 - x) (1 - x^2) (1 - x^3) (1-x^4) ... \}^6} , 
\\[1ex]
p(5) + p(12) x + p (19)x^2 + p(26) x^3 + p (33) x^4 + ... 
\\
\qquad =7 {\{(1-x^7) (1-x^{14}) (1-x^{21}) (1 - x^{28}) ... \}^3 \over 
\{(1-x) (1-x^2) (1 -x^2) (1-x^3) (1 - x^4) ... \}^4 } 
\\
\qquad\quad {}+49 {x\{(1-x^7) (1 - x^{14}) (1 - x^{21}) (1 - x^{28}) ... \}^7 \over
\{(1-x) (1-x^2)(1-x^3) (1-x^4) ... \}^8}, 
\end{gather*}
is worth paying attention to.
 
What other interesting things has the magician discovered?  Well, we can begin with the following 
quotations:``
\begin{gather*}
\int^\infty_{0} { \sin(nx)dx \over {\dps{x + {1 \over x} { \atop +}  {2 \over x} { \atop +} 
{3 \over x} { \atop +}  { \atop ...} }}} = {\sqrt({1 \over 2} \pi) \over {\dps {n  + 
{1 \over n} { \atop +}  {2 \over n} { \atop +}  {3 \over n} { \atop +} { \atop ...} } }} ;
\\[1ex]
\int^\infty_{0} {\sin({1 \over 2} \pi n x)  d x \over {\dps{x + {1^2 \over x} { \atop +} 
{2^2 \over x} { \atop +} {3^2 \over x} { \atop +} { \atop ...}}}}  = {1 \over n} { \atop +} 
{1^2 \over n} { \atop +} {2^2 \over n} { \atop +}  {3^2 \over n} { \atop +} { \atop ...} ; 
\\[1ex]
4 \int^\infty_{0} {xe^{-x\sqrt{5}} \over \cosh x} dx = {1 \over 1} { \atop + } {1^2 \over  1} 
{ \atop +} {1^2 \over 1} { \atop +}  {2^2 \over 1} { \atop +}  {2^2 \over 1} { \atop +} 
{3^2 \over 1} { \atop +}  {3^2 \over 1} { \atop +} { \atop ...} ; 
\\[1ex]
2 \int^\infty_{0} {x^{2} e^{-x\sqrt{3}} \over \sinh x} dx = {1 \over 1} { \atop +}  {1^3 \over 1} 
{ \atop +} {1^3 \over 3} { \atop +} {2^3 \over 1} { \atop +} {2^3 \over 5} { \atop +} 
{3^3 \over 1} { \atop +}  {3^3 \over 7} { \atop +} { \atop ...} ; 
\end{gather*}

If ${1 \over 2} \pi \alpha = \log\tan\{{1 \over 4} \pi (1 + \beta)\}$, then
\begin{equation*}
{1^2 + \alpha^2 \over 1^2 - \beta^2} \bigg({3^2 - \beta^2 \over 3^2 + \alpha^2} \bigg)^3 \bigg({ 
5^2 + \alpha^2 \over 5^2 - \beta^2} \bigg)^5 \bigg({7^2 - \beta^2 \over 7^2 + \alpha^2} \bigg)^7
\bigg({9^2 + \alpha^2 \over 9^2 - \beta^2}\bigg)^9...= e^{\pi \alpha \beta/2}; 
\end{equation*}
 
If 
\begin{equation*}
F (x) = {1 \over 1} { \atop +} {x \over 1} { \atop +} {x^2 \over 1} { \atop +} {x^3 \over 1} 
{ \atop +} {x^4 \over 1} { \atop +} {x^5 \over 1} { \atop +} { \atop ...} , 
\end{equation*}
then
\begin{equation*}
\bigg\{ {\sqrt{5} + 1 \over 2} + e^{-2\alpha/5} F(e^{-2\alpha} ) \bigg\} \bigg\{{\sqrt{5 }+1 
\over 2} + e^{-2 \beta/5} F(e^{-2 \beta} ) \bigg\} = {5 + \sqrt{5} \over 2} 
\end{equation*}
with the condition $\alpha \beta = \pi^2$. \\
The above theorem is a particular case of a theorem on the continued fraction
\begin{equation*}
{1 \over 1} { \atop +} {ax \over 1} { \atop +} {ax^2 \over 1} { \atop +}
 {ax^3 \over 1} { \atop +} {ax^4 \over 1} { \atop +} {ax^5 \over 1} { \atop +} { \atop ...} , 
\end{equation*}
which is a particular case of the continued fraction
\begin{equation*}
{1 \over 1} { \atop +} {ax \over 1 + bx} { \atop +} {ax^2 \over 1 + bx^2} { \atop +}
 {ax^2 \over 1 + bx^3} { \atop +} { \atop  ...} , 
\end{equation*}
which is a particular case of a general theorem on continued fractions;
\begin{gather*}
{a \over 1+n} { \atop +} {a^2 \over 3+n} { \atop +} {(2a)^2 \over 5 + n} { \atop +}
 {(3a)^2 \over 7 + n} { \atop +} { \atop ...} 
\\
\qquad = 2a \int^1_{0}  z^{n/ \sqrt{1+a^2}} {dz \over \big({\sqrt{ 1+ a^2} + 1 \big)+z^2 \big( 
\sqrt{1 + a^2 }- 1 \big) }} ,
\end{gather*}
which is a particular case of the continued fraction
\begin{equation*}
{a \over p + n} { \atop +} {1 \cdot p \cdot a^2 \over p + n + 2} { \atop +} {2 (p + 1)a^2 \over 
p + n + 3} { \atop +} { \atop ...} , 
\end{equation*}
which is a particular case of a corollary to a theorem on transformation of integrals and 
continued fractions.''

The formulae above are taken from Ramanujan's {\it{Collected Papers}} [26].  
Now the story gets complicated.  From his schooldays, Ramanujan had been recording his 
mathematical discoveries in notebooks he kept through the rest of his short life.  He was 
too poor to afford writing paper;  so he worked on a slate and then copied the final
{\it{results}} $-$ without any {\it{proofs}} $-$ into the notebooks; year after year, from about 
1903, when he was 16, until his death in 1920.

What was to be done with the notebooks?  Those who have seen them, Hardy and Littlewood 
in particular, have recognized the notebooks for the treasure trove they were and urged their 
edition and publication.  In 1929, G. N. Watson and B. M. Wilson undertook the editing task, 
estimating it to take five years to complete.  This was not to be.  Wilson died in 1935 after  
routine surgery, and Watson eventually abandoned the project in the late 1930's.  On the 
Indian side, things were scarcely  better.  Ranganathan [28] relates, matter of factly, 
the shameful indifference which met his decades-long attempts to have Ramanujan's 
notebooks published.  Finally, in 1957, an unedited photostat copy was produced by the Tata
 Institute in Bombay; but anyone who has seen this copy can testify that it is very hard to make 
sense out of and it is near useless for all 
but dedicated Ramanujan scholars.  

Here matters stood; and stood; for a long-time; awaiting their next chance at the Ramanujan 
birthday centennial approaching in 1987; while the civilized word meanwhile indulged in the 
annual orgy celebrating the birthdays of their biggest mass murderers aka greatest national 
heroes:  Lincoln in the USA, Lenin in the Soviet Union, ...

Later rather than sooner, luck finally intervened.  George E. Andrews in 1976 
 discovered Ramanujan's ``lost'' notebook, recognized it for the treasure it was, and published 
some parts of it [1-3].  This has evidently rekindled the semi-dormant interest 
in Ramanujan, kept alive by a small circle of practicing partitioners, 
 into something like a post-critical flaming explosion of attention, 
and the rest is history $-$ albeit very recent 
one, and the subject of this review.

Bruce C. Berndt has been continuously working through the Ramanujan notebooks since May 1977, 
having devoted all of his research time to proving the claims made by Ramanujan in his 
notebooks.  By Berndt's count, the total number of results examined in five volumes is an astounding 
3254.  Is there anyone among us who would not be daunted by the enormity of such a task?

(Berndt: ``In notes left by B. M. Wilson, he tells us how George Pol$\acute{\rm{y}}$a 
was captivated by Ramanujan's 
formulas.  One day in 1925 while Pol$\acute{\rm{y}}$a  was visiting Oxford, he borrowed 
 from Hardy his copy 
of Ramanujan's notebooks.  A couple of days later, Pol$\acute{\rm{y}}$a 
 returned them in almost a state of 
panic explaining that however long he kept them, he would have to keep attempting to verify 
formulas therein and never again would have time to establish another original result of 
his own.'')

How difficult was the task?  G. N. Watson delivered a lecture entitled   
{\it{Ramanujan's Note Books}} at the February 5, 1931,  meeting of the London Mathematical 
Society [29].  
Taking from the notebooks (p. 150, [29]) a pair of modular equations of order 13, 
he went on to say: ``This pair of formulae took me a month to prove, though 
I am now fairly certain that my proof 
is the same as Ramanujan's, with the exception of one section, which I think 
that he was able to work 
out more neatly than I have succeeded in doing.  He has given a somewhat simular form of the 
modular equation of order seventeen which I have not yet worked out, though it will probably 
prove amenable to the same methods; and a rather more complex form for the modular equation 
of order nineteen about which I am less hopeful.''

Complicating matters, one had to verify carefully {\it{every}} claim made by Ramanujan, for some of these have 
misprints or simple mistakes, while others are seriously wrong; and the  handful of the latter 
have so far resisted all attempts to be salvaged. 

There is an old dictum, ``Every man is a debtor to his 
profession.''  Such a debt is the most obvious in mathematics 
where every past discovery provides a 
perpetual free lunch for future generations, until the end of time or civilization $-$ 
whichever comes first. 
It is a tribute to Berndt's talents and perseverance that he has completed this rather 
terrifying project.  He is a credit to mathematics and a true benefactor of mankind.

Which of the Ramanujan formulae are the most interesting?  The question is clearly meaningless, 
for the answer will vary greatly from one reader to another.  Less obviously, the answer 
will likely change in time for a fixed reader.  In preparing this review, I have re-read  all 
five of Brendt's volumes, and observed that about half of the entries I had marked $-$ years ago on 
first readings, as deserving further study $-$ no longer seem as interesting; while an 
 approximately 
similar amount of entries previously glossed over now show clear promise of profit when 
analyzed.  One's interests, expertise, and time available change; but Ramanujan's formulae 
endure forever.

Below are listed a few formulae $-$ out of thousands $-$ from the Notebooks.  They have been 
chosen as likely to be fascinating to a general reader. 
To counter-balance my inevitably subjective choice, I didn't include any formula that seemed 
promising to my personal research interests.  So:

\begin{gather*}
\int^\infty_0 {dx \over (x^2 + 11^2) (x^2 + 21^2) (x^3 + 31^2) (x^2 + 41^2) (x^2 + 51^2)}  
\\
\qquad = {5 \pi \over 12 \cdot 13 \cdot 16 \cdot 17 \cdot 18 \cdot 22 \cdot 23 \cdot 24 \cdot 
31 \cdot 32 \cdot 41} \ ; 
\\[1ex]
\pi = 3.14159265358... 
\\[1ex]
\bigg(9^2 + {19^2 \over 22} \bigg)^{{1 \over 4}} = 3.14159265262...
\\[1ex]
{355 \over 113} \bigg(1 - {.0003 \over 3533} \bigg) = 3.141592653589743... \ ; 
\end{gather*} 
 
For $n > 0$,
\begin{equation*}
\sum^\infty_{k=1} \tan^{-1} \bigg({2 \over (n + k)^2} \bigg) = \tan^{-1} \bigg({2n+1 \over 
n^2 + n-1} \bigg) + \rho (n), 
\end{equation*}
where $\rho(n) = \pi$ if $n < (\sqrt{5} -1)/2$ and $\rho(n) = 0$ otherwise;
\begin{equation*}
\sum^\infty_{k=0} (-1)^k {(2k + 1)^3 + (2k + 1)^2 \over k!} = 0; 
\end{equation*}
 
For $|x|<1$, 
\begin{equation*}
\Pi^\infty_{k=1} (1 - x^{p_{k}})^{-1} = 1 + \sum^\infty_{k=1} {x^{p_1 + p_2 + ... + p_k} \over 
(1-x) (1-x^2) ... (1 - x^k)}, \tag{*}
\end{equation*}
where $p_1, p_2$,... denote the primes in ascending order. \\
This entry, in fact, is canceled by Ramanujan.  Let $c_n$ and $d_n, \ 2 \leq n < \infty,$ denote the 
coefficients of $x^n$ on the left and right sides, respectively, of \ $(*)$.  Then, quite amazingly, 
$c_n = d_n$ for $2 \leq n \leq 20$.  But $c_{21} = 30$ and $d_{21} = 31$.  Thus, as indicated by 
Ramanujan, \ $(*)$ is false;
\begin{equation*}
e^{n\sin^{-1}x} = 1 + nx + \sum^\infty_{k=2} {b_k (n)x^k \over k!}, \ \ |x| \leq 1, 
\end{equation*}
where, for $k \geq 2$, 
\begin{gather*}
b_k (n) = \bigg\{ 
\matrix{%
n^2 (n^2 + 2^2) (n^2 + 4^2)&... &(n^2 + (k - 2)^2), &{\rm{if}} \ k \ 
{\rm{is \ even}}, \cr
n(n^2 + 1^2) (n^2 + 3^2) &... &(n^2 + (k - 2)^2), & {\rm{if}} \ k \ {\rm{is \ odd;}}\cr}
\\[1ex]
e^{ax} = 1 + {a \sin (cx) \over c} + \sum^\infty_{k=2} {u_k \over k!} \bigg({\sin(cx) \over c} 
\bigg)^k, 
\end{gather*}
where, for $k \geq 2$, 
\begin{equation*}
u_k = \bigg\{\matrix{a^2 \{a^2 + (2c)^2\} \{a^2 + (4c)^2\}&...&\{a^2 + (k-2)^2 c^2\},&  
{\rm{if}} \ k \ {\rm{is \ even}}, \cr 
a\{a^2 + c^2\} \{a^2 + (3c)^2\}&... &\{a^2 + (k - 2)^2c^2\}, &{\rm{if}} \ k \ 
{\rm{is \ odd;}}\cr} 
\end{equation*}
 
For a given function $F(x)$, denote $F^0 (x) = x, \ F^1 (x) = F(x), \ F^{n+1} (x) = 
F(F^n(x)).$  Extend this definition into the family $F^r(x), $ for all $r \in {\bf{R}}$.

Let $F(x) = x^2 - 2, \ x \geq 2$.  Then  
\begin{gather*}
F^{1/2} (x) = \bigg({x+ \sqrt{x^2 - 4 \over 2}} \bigg)^{\sqrt{2}} + 
\bigg({x-\sqrt x^2 - 4 \over 
2 } \bigg)^{\sqrt{2}}, 
\\[1ex]
F^{\Log 3/\Log 2} (x) = x^2 - 3x, 
\\[1ex]
F^{\Log 5/\Log 2} (x) = x^5 - 5x^3 + 5x; 
\end{gather*}

Let $0 < x \leq 1$.  If
\begin{equation*}
1 + \sqrt{F^{\Log 2} (x)}  = \sqrt{{1 - F^{\Log 2}(x) \over 1-x}}, 
\end{equation*}
then
\begin{equation*}
1 + 2 \sqrt{{F^{\Log 3}(x) \over x}} = \sqrt{{1-F^{\Log 3}(x) \over 1-x}} ;
\end{equation*}
 
If $Re \ x > {1\over 2}$, then
\begin{equation*}
1 + 3 {x-1 \over x+1} + 5 {(x-1) (x-2) \over (x+1) (x+2)} + ... = x; 
\end{equation*}
 
If $Re \ x> 1$, then
\begin{equation*}
1 - 3 {x-1 \over x+1} + 5 {(x-1) (x-2) \over (x+1) (x+2)} + ... = {x \over 
2x - 1} ; 
\end{equation*}
 
If $Re \ x > {3\over 2}$, then
\begin{gather*}
1^3 + 3^3 {x-1 \over x+1} + 5^3 {(x-1) (x-2) \over (x+1) (x+2)} + ... =x (4x-3); 
\\[1ex]
\bigg\{1 + {1^2 + n \over 4^2}  x + {(1^2 + n) (5^2 + n) \over 4^2 \cdot 8^2} x^2 + ... 
\bigg\}^2 
\\*
\qquad = 1 + {1 \over 2} {1^2 + n \over 2^2} x + {1 \cdot 3 \over 2 \cdot 4} {(1^2 + n) (3^2 + n) 
\over 2^2 \cdot 4^2} x^2 + ... ; 
\\[1ex]
3 = (1 + 2 (1 + 3 (1 + 4 (1 + ...) ^{1/2})^{1/2})^{1/2})^{1/2} ; 
\\[1ex]
2 = (6 + (6 + (6 + ...)^{1/3})^{1/3})^{1/3};
\\[1ex]
{5 \over 3} = {4 \over 1} { \atop +}  {6 \over 3} { \atop +} 
 {8 \over 5} { \atop +}  {10 \over 7} { \atop +}  {\atop ...} ;
\end{gather*}
 
As a formal identity,
\begin{equation*}
{a_1 + h \over 1} { \atop +} {a_1 \over x} { \atop +}  {a_2 + h \over 1} { \atop +} 
{a_2 \over x} { \atop + ...}  = h + {a_1 \over 1} { \atop +} {a_1 +h \over x} { \atop +} 
{a_2 \over 1} { \atop +}  {a_2 + h \over x} { \atop + ...} ; 
\end{equation*}
 
We have
\begin{equation*}
lim_{x \rightarrow + \infty} \bigg\{ \sqrt{{2x \over \pi}} - {x \over 1} { \atop +} 
{2x \over 2} { \atop +}  {3x \over 3} { \atop +}  {4x \over 4} { \atop + ...}  \bigg\} = 
{2 \over 3 \pi}; 
\end{equation*}
 
If $Re \ x > 0$, then
\begin{equation*}
2 \sum^\infty_{k = 0} {1 \over (x + 2k + 1)^2} = {1 \over x} { \atop +}  {1^4 \over 3x} 
{ \atop +} {2^4 \over 5x} { \atop +} {3^4 \over 7x} { \atop + ...} ; 
\end{equation*}
 
If $Re \ x > 0,$ then
\begin{equation*}
2 \sum^\infty_{k=0} {(-1)^k \over (x + 2k + 1)^2} = {1 \over x^2 -1} { \atop +} 
{2^2 \over 1} { \atop +}  {2^2 \over x^2-1} { \atop +} {4^2 \over 1} { \atop +}  {4^2 \over 
x^2 - 1} { \atop + ...} ; 
\end{equation*}
 
If $n$ is any complex number outside of $(- \infty,0]$, then
\begin{gather*}
\int^\infty_{0} e^{-x} (1 + x/n)^n dx 
\\
\qquad = 1 + {n \over 1} { \atop +}  {1 (n-1) \over 3} { \atop +}  {2(n-2) \over 5} { \atop +} 
{3(n-3) \over 7} { \atop + ...}  
\\
\qquad = 2 + {n-1 \over 2} { \atop +} {1(n-2) \over 4} { \atop +}  {2(n-3) \over 6} { \atop +} 
{3(n-4) \over 8} { \atop + ...} ; 
\end{gather*}
 
Let $\alpha, \beta > 0$, with $\alpha \beta = \pi^2$.  Then
\begin{equation*}
e^{(\alpha-\beta)/{1 2}} = \bigg({\alpha \over \beta}\bigg)^{1/4} \Pi^\infty_{k=1} 
{1 - e^{2\alpha k} \over 1-e^{2\beta k}} ;
\end{equation*}
 
If $|q| < 1,$ then
\begin{gather*}
\sum^\infty_{k=0} (1 - a)^k q^{k(k+1)/2} = {1\over 1} { \atop +} {aq \over 1} { \atop +} 
{a(q^2 - q) \over 1} { \atop +}  {aq^3 \over 1} { \atop +}  {a(q^4 - q^2) \over 1} 
{ \atop + ...} ;
\\[1ex]
{11 \over 10} {1111 \over 1110} {111111 \over 111110} ... = 1.101001000100001...;
\\[1ex]
{e^{-\pi/5} \over 1} { \atop -}   {e^{-\pi} \over 1} { \atop +} {e^{-2\pi} \over 1} { \atop - ...} 
 = \sqrt{{5 - {\sqrt{5}} \over 2}} - {\sqrt{5} - 1 \over 2} ;
\\[1ex]
{e^{-2\pi/5} \over 1} { \atop +} {e^{-2 \pi} \over 1} { \atop +}  {e^{-4 \pi} \over 1} 
{ \atop + ...}  = {\sqrt{5+ \sqrt{5}} \over 2} - {\sqrt{5} +1 \over 2};
\\[1ex]
\sum^\infty_{k=1} (k^2 \pi - {1 \over 4}) e^{-\pi k^{2}} = {1 \over 8}; 
\end{gather*}
 
For $n > 0$, Re $\alpha >0$, and Re $\beta > 0$, define
\begin{equation*}
F(\alpha, \beta) = {\alpha \over n} { \atop +} {\beta^2 \over n} { \atop +}  {(2\alpha)^2 
\over n} { \atop +}  {(3 \beta)^2 \over n} { \atop +}  {(4 \alpha)^2 \over n} { \atop + ...} 
. 
\end{equation*}
Then
\begin{equation*}
F\bigg({\alpha + \beta \over 2}, \sqrt{\alpha \beta} \bigg) = {1 \over 2} (F(\alpha, \beta) + 
F (\beta, \alpha)); 
\end{equation*}
 
If we approximate $\pi$ by (97$^{1/2}-{1 \over 11})^{1/4}$ in the expression 
${1\over 2}d \sqrt{\pi}$, then if a circle of one million square miles is taken, the error 
made is approximately 1/100th of an inch;
 
If $n > 0$, then
\begin{equation*}
\int^n_{0} \bigg({n \over x} \bigg)^x dx = \sum^\infty_{k = 1} {n^k \over k^k} ;
\end{equation*}
 
Let $\varphi (x)$ be continuous on $[0, \infty)$ and suppose that $lim_{x \rightarrow \infty} 
\varphi (x) = : \varphi (\infty)$ exists.  Then
\begin{equation*}
\lim_{n\rightarrow 0} n \int^\infty_{0} x^{n-1} \varphi (x) dx = \varphi (0) - 
\varphi (\infty); 
\end{equation*}
 
We have
\begin{gather*}
{1 \over 4} + 2 = (1 {1 \over 2})^2, 
\\
{1 \over 4} + 2 \cdot 3 = (2 {1\over 2})^2, 
\\
{1\over 4} + 2 \cdot 3 \cdot 5 = (5 {1 \over 2})^2, 
\\
{1\over 4} + 2 \cdot 3 \cdot 5 \cdot 7 = (14 {1\over 2})^2, 
\\
{1\over 4} + 2 \cdot 3 \cdot 5 \cdot 7 \cdot 11 \cdot 13 \cdot 17 
 = (714{1\over 4} )^2; 
\end{gather*}
 
The expressions
\begin{gather*}
{12 \over {\sqrt{130}}} \log \bigg({(3 + \sqrt{13)} (\sqrt{8} + \sqrt{10}) \over 2} \bigg), 
\\[1ex]
{24 \over {\sqrt{142}}} \log \bigg({\sqrt{10 + 11 \sqrt{2}} + {\sqrt{10 + 7 \sqrt{2}}} \over 2} 
\bigg), 
\end{gather*}
and
\begin{equation*}
{12 \over \sqrt{190}} \log ((3 + \sqrt{10}) (\sqrt{8} + {\sqrt{10}})) 
\end{equation*}
are approximates to $\pi$ valid for $14, 15,$ and 18 decimal places, respectively.

And one entry from the ``lost'' notebook [27, 18]:
 
If
\begin{gather*}
\sum_{n \geq 0} a_n x^n = {1 + 53 x + 9 x^2 \over 1 - 82 x - 82 x^2 + x^3}, 
\\[1ex]
\sum_{n \geq 0} b_n x^n = {2 - 26 x - 12 x^2 \over 1 - 82 x - 82 x^2 + x^3}, 
\end{gather*}
and
\begin{equation*}
\sum_{n \geq 0} c_n x^n = {2 + 8x - 10 x^2 \over 1 - 82 x - 82 x^2 + x^3}, 
\end{equation*}
then
\begin{equation*}
a^3_{n} + b^3_n = c^3_n + (-1)^n. 
\end{equation*}

A few comments are in order.

In 1929, reviewing Ramanujan's {\it{Collected Papers}}, Littlewood remarked [23]:  ``But the 
great day of formulae seems to be over.  No one, if we are again to take the highest  
standpoint, seems able to discover a radically new type, though Ramanujan comes near it in his 
work on partition series..."  The great Littlewood seems to be unduly pessimistic there, as the 
days of interesting new formulae are still with us, and the days of {\it{great}} formulae 
seem to have hardly began.  

(In his Commonplace Book [14, p. 61], Littlewood relates:  ``I read in the proof-sheets of 
Hardy on Ramanujan:  `As someone said, each of the positive integers was one of his personal 
friends.' My reaction was, `I wonder who said that; I wish I had.'  In the next proofsheets I 
read (what now stands), `It was Littlewood who said...' \ '')

There are many mysteries surrounding Ramanujan, some of which may remain un-understood forever, 
such as his {\it{methods}}.  Berndt's heroic work has been devoted to {\it{verification}} 
of Ramanujan's results, and there are quite a number of these where it remains totally opaque 
as to 
just what the route Ramanujan used to arrive at his beautiful formulae.  There are also minor 
mysteries, of less import, such as why Ramanujan got no substantial official help or 
recognition in India in the years up until shortly before 
his genius was recognized publicly by Hardy.  The accepted 
line, that there was no one in India at the time competent enough to understand any of his 
work, seems to me untenable, for the listing of questions Ramanujan had submitted to the Problems 
section of 
{\it{The Journal of Indian Mathematical Society}} (printed on pp. 323-334 of his {\it{Collected 
Papers)}} shows that almost all of the problems he had proposed have been solved by a variety of 
people.

Ramanujan's talents thus {\it{were}} appreciated in India, just not by anyone who was in a 
position to  help him, for people in power, then as now, have more important problems to 
occupy them, the fate of the Universe being their major concern.

It would have been very improbable for no one in Ramanujan's surroundings being able to comprehend 
his gifts, and indeed he did have a small handful of devoted friends and admirers who had 
managed to prevent Ramanujan from starving to death $-$ but just barely.  He was similarly 
unfortunate in the choice of his parents, but not singularly so, for dim people are always 
and everywhere in the majority.  At Ramanujan's funeral on April 26, 1920, most of his orthodox 
Brahmin relatives stayed away, since Ramanujan was tainted in their eyes, having crossed  the 
waters on his sea voyage to England in April 1914.

This last tidbit comes courtesy of R. Kanigel's biography of Ramanujan [19].  There 
are other interesting data in that book, but it is ultimately unsatisfying $-$ and I don't 
mean the modern abominations of biography writing, such as chercherizing les femmes (or 
les hommes in this particular case); the fault lies deeper, for one is interested in 
Ramanujan only because of his {\it{mathematics}} and nothing else, so a proper biography of 
Ramanujan, were it to be written, should be a {\it{mathematical biography}}.  Meanwhile, R. 
Askey's insightful review [7] could serve as interesting ersatz, chronicling in 
detail how Ramanujan had rapidly rediscovered most of the results of every classical subject 
he put his mind to.  After a few weeks of perusing Ramanujan's {\it{Collected Papers}} and 5-volume {\it{Notebooks}}, 
one gets a distinct feeling that Ramanujan's real talent was not, as Hardy and others thought, his 
tremendous insight into algebraic transformations and so forth, but the true mathematical gift 
of multidimensional perception, of gaze penetrating behind the superficial surface of things, of 
undogmatic rejection of accepted meanings and conventions.

(Berndt: ``Paul Erdos has passed on to us Hardy's  personal ratings of mathematicians.  Suppose 
that we rate mathematicians on the basis of pure talent on a scale from 0 to 100.  Hardy 
gave himself a score of 25, Littlewood 30, Hilbert 80, and Ramanujan 100.'')

The adventure of the human existence can be justified by the rare instances of inspired 
human conduct and by the glory of the inventive human mind.  The latter through the centuries has
given us the {\it{Magna Carta}}, the {\it{Declaration of Independence}}, the {\it{Bill of 
Rights}}, and, with the trivial changes of sign, marxism, communism, socialism, and all the 
other perversims.  Ramanujan's discoveries, now made available to all through Berndt's 
decades-long endeavor, 
lift us another notch on the eternal climb in search of the Ultimate Truth.

To conclude this part of the review of Brendt's monumental 5-volume set, perhaps I could do no 
better than to quote Auberon Waugh's advice to the readers of his book {\it{Way Of The World}}:
  ``I should warn that the density and richness of the material make it unsuitable for 
prolonged reading.  Perhaps the best plan is to have a copy in every room for immediate 
relief in moments of desolation, as well as several at your place of work to impress colleagues 
and visitors.''

\subsection*{From Ramanujan to Fermat}

\begin{flushright}
\parbox{5.5cm}{\sffamily \small
Things are seldom what they seem; 
Skim milk masquerades as cream... \\
\null\qquad  Gilbert \& Sullivan}
\end{flushright}
\smallskip
Hardy observed:  ``There is always more in one of Ramanujan's formulae than meets the eye, as anyone who sets 
to work to verify those which look the easiest will soon discover.  In some the interest lies 
very deep, in others comparatively near the surface; but there is not one which is not 
curious and entertaining.''

The last entry in Chapter 4 of Ramanujan's 2$^{nd}$ Notebook [8, p.\ 108], reads: 
\begin{entry}{15}
For each positive integer $\kappa$, let $G_k = \sum_{0 \leq 2n+1 \leq k} 
1/(2n+1)$.  Then for all complex $x$,
\begin{equation*}
e^x \sum^\infty_{k=1} {(-2)^{k-1} x^k \over k!k} = \sum^\infty_{k=1} 
{G_k x^k \over k!} . \tag{1}
\end{equation*}
\end{entry}

This identity is certainly not deep, and on first sight it doesn't look curious or entertaining 
either.  But with Ramanujan one never knows for sure.  Berndt provides a very short proof of 
(1) and then, as if in an afterthought, remarks:  ``Note that by equating coefficients of $x^n$  
in (1), we find that
\begin{equation*}
\sum^n_{k=1} \bigg({n \atop k}\bigg) {(-2)^{k-1} \over k} = G_n, \ \ \ n \geq 1. {\rm{"}}
 \tag{2}
\end{equation*}

Ahh, but this equality {\it{is}} interesting:  a regular binomial sum on the left represents  
 a {\it{step-function}} on the right which doesn't change 
when $n$ increases by 1 whenever $n$ is even.  This is reminiscent of some formulae in 
$q$-series (see, e.g. [22]) where a different form of ``2-periodicity'' shows up, such as in 
the Gauss formula
\begin{equation*}
\sum^N_{\ell = 0} \bigg[ {N \atop \ell} \bigg]_q (-1)^\ell =
\bigg\{\matrix{0, \ \ \  N \ {\rm{odd}} > 0, \cr
(1 - q) (1-q^3) ... (1 - q^{2m-1} ), \ \ N = 2m > 0; \cr}
\end{equation*}
here $\dps{\bigg[{N \atop \ell}\bigg] = \bigg[{N \atop \ell}\bigg]_q} $ is the $q$-binomial 
coefficient:
\begin{equation*}
\bigg[{\gamma \atop \ell}\bigg] = {[\gamma] ... [\gamma - \ell + 1] \over [\ell] ... 
[1] },  \ \ \ell \in \N; \ \ \ \bigg[{\gamma \atop 0}\bigg] = 1: \ \ [\gamma] = [\gamma]_q 
= {1 - q \over 1 - q}^\gamma. \tag{3}
\end{equation*}

It could be instructive to find a (``quantum'') $q$-analog of the equality (2), and this is 
what we shall now attempt.  How should we proceed?  There are no set rules for $q$-tization, 
for an infrastructure of $q$-mathematics is lacking, and it may not even exist in the 
classical sense.  Certainly, trying to $q$-deform separately each entry in a given equality 
to  be quantised is a fool's errand, tedious in the extreme; and in our particular case (2) 
this errand can not even be attempted with confidence thanks to the presence of the factor 
$2^{k-1}$ in the LHS:  in contrast to algebraic {\it{expressions}} such as 
$\dps{{n \choose k}}$, the fixed {\it{numbers}} such as 2 turn under quantization 
into unpredictable situation-dependent objects, frequently sprinkled with 
classically-invisible factors ${\dps{{1+a \over 1+b}}}$, where $a$ and $b$ are polynomials in $q$ attaining 
equal values for $q=1$.  A more intelligent approach would be to notice that binomial 
coefficients dominate our equality (2), rewrite it in hypergeometric notation, and then convert 
 hypergeometric objects into basic $(=q-)$ 
hypergeometric ones;  this approach to quantization was advocated by Andrews [4], and 
it works well sometimes  [4] and not so well at other times [5]; it couldn't be attempted 
in our case (2), again because of the $2^{k-1}$-factor.

There is nothing then left to do but revert to the remedy of last resort:  to try to quantize the 
{\it{argument}} used to establish the equality we are trying to quantize.  Berndt gives the 
following proof of formula (1):

We have
\begin{equation*}
I \equiv e^x \int^1_0 {1 - e^{-2xz} \over 2z}  dz = e^x \int^1_0 \sum^\infty_{k=1} 
{x^k (-2z)^{k-1} \over k!} dz = e^x \sum^\infty_{k=1} {(-2)^{k-1} x^k \over k!k} . \tag{4a}
\end{equation*}
On the other  hand,
\begin{align*}
I &= \int^1_0 {e^x-e^{x(1-2z)} \over 2z} dz = \sum^\infty_{k=1} {x^k \over k!} \int^1_0 
{1-(1-2z)^k \over 2z} dz 
\\
&= {1 \over 2} \sum^\infty_{k=1} {x^k \over k!} \int^1_{-1} {1 - t^k \over 1-t} dt = 
\sum^\infty_{k=1} {x^k \over k!} \int^1_0  \sum_{0 \leq j \leq (k-1)/2} t^{2j} dt = \sum^\infty
_{k=1} {G_kx^k \over k!} . \tag{4b}
\end{align*}
Combining (4a) and (4b), we deduce (1).

One of the steps involved in the derivation above, and the only one where the direct quantization 
breaks down, is the equality
\begin{equation*}
e^x e^{-2xz} = e^{x(1-2z)} . \tag{5}
\end{equation*}
$q$-exponentials have many nice properties, but
\begin{equation*}
e^x e^y = e^{x+y}, \ \  \ \ \ xy = yx,  
\end{equation*}
is not one of them.  What to do now?  Well, these pesky exponentials appear inevitably 
because Ramanujan has chosen for his formula (1) to work with the {\it{exponential}} 
generating function for the sequence $\{G_n\}$; since we are after not the identity (1) 
itself, but only the equality (2), we should find a proof of this equality based on 
the straightforward generating function for the $G_n$'s.  So:
\begin{gather*}
\sum^\infty_{N=1} G_N x^N = \sum^\infty_{m=0} x^{2m+1} \sum^m_{s=0} {1 \over 2s+1} 
+ \sum^\infty_{m=1} x^{2m} \sum^{m-1}_{s=0} {1 \over 2s+1} \tag{6a}
\\[1ex]
\qquad = \sum^\infty_{s=0} {1 \over 2s+1} \sum^\infty_{k=0} \bigg(x^{2s+1+2k} + x^{2s+2+2k} \bigg) 
= \sum^\infty_{s=0} {x^{2s+1} \over 2s+1} (1 + x) \sum^\infty_{k=0} x^{2k} \tag{6b}
\\[1ex]
\qquad= \sum^\infty_{s=0} {x^{2s+1} \over 2s+1} (1 + x) {1 \over 1-x^2} = {1 \over 1-x} 
\sum^\infty_{s=0} {x^{2s+1} \over 2s+1} 
\\
\qquad= {1 \over 1-x} {\log (1+x) - \log(1-x) \over 2} 
\tag{6c}
\\[1ex]
\qquad={1 \over 2(1-x)} \log \bigg({1 +x \over 1-x}\bigg) = {1 \over 2(1-x)} \log \bigg(1 + {2x \over 
1-x}\bigg) 
\\
\qquad= {1 \over 2(1-x)} \sum^\infty_{k=1} \bigg({2x \over 1-x}\bigg)^k {(-1)^{k-1} \over 
k} \tag{6d}
\\[1ex]
\qquad= \sum^\infty_{k=1} {x^k \over k} {(-2)^{k-1} \over (1 - x)^{k+1}} = \sum^\infty_{k=1} 
{x^k \over k} (-2)^{k-1} \sum^\infty_{n=0} \bigg({k+n \atop k}\bigg) x^n 
\\
\qquad= \sum^\infty_{N=1} 
x^N \sum^N_{k=1} {(-2)^{k-1} \over k} \bigg({N \atop k}\bigg). \tag{6e}
\end{gather*}
Comparing the first and the last entries in this chain of identifies, we recover formula (2).  
Notice, that on the way we have also found that
\begin{equation*}
\sum^\infty_{N=1} G_N x^N = {1 \over 2(1-x)} \log \bigg({1 + x \over 1-x} \bigg). \tag{7}
\end{equation*}

Looking over each of the 5 lines $(6a-e)$ above, we see that only one of them, $(6d)$, throws 
up some obstructions to quantization, namely the trio of equalities

\begin{equation*}
\log(1+x) - \log (1-x) = \log \bigg({1 + x \over 1-x} \bigg) = \log \bigg(1 + {2x \over 1-x}
\bigg) 
= \sum^\infty_{k=1} \bigg({ 2x \over 1-x}\bigg)^k {(-1)^{k-1}  \over k} . \tag{8}
\end{equation*}
If we agree to ignore the two intermediate identities in the chain (8), we need only to find a 
$q$-analog of the equality
\begin{equation*}
\log (1-x) - \log (1-x) = \sum^\infty_{k=1} \bigg({2x \over 1-x}\bigg)^k {(-1)^{k-1} \over k}; 
\tag{9}
\end{equation*}
equivalently, we need a $q$-analog of the equality
\begin{equation*}
\sum^\infty_{s=0} {x^{2s+1} \over 2s+1} = \sum^\infty_{k=1} \bigg({x \over 1-x}\bigg)^k 
{(-2)^{k-1} \over k} ; \tag{10}
\end{equation*} 
differentiating this with respect to $x$ and using the obvious relation
\begin{equation*}
{d \over dx} \bigg(\bigg({x \over 1-\alpha x} \bigg)^k \bigg) = k {x^{k-1} \over (1-\alpha x)^{k+1}}, 
\tag{11}
\end{equation*}
we get still another form of the relation (9):
\begin{equation*}
\sum^\infty_{s=0} x^{2s} = \sum^\infty_{k=0} {(-2x)^k \over (1-x)^{k+2} }. \tag{12}
\end{equation*}
\begin{proposition}{13}  
We have
\begin{equation*}
\sum^\infty_{s=0} x^{2s} = \sum^\infty_{k=0} {(-1)^k (1 \dot + q)^k x^k \over 
(1 -^{\hskip-.09truein \cdot} x)^{k+2}}, \tag{14}
\end{equation*}
where
\begin{equation*}
(u \dot + v)^k = \Pi^{k-1}_{i=0} (u + q^i v), \ \ k \in \N; \ \ (u \dot + v)^0 = 1. 
\tag{15}
\end{equation*}
\end{proposition}
\begin{proof}
Multiplying both sides of (14) by $(1-x)$, 
we reduce formula (14) to the identity 
\begin{equation*}
{1 \over 1+x} = \sum^\infty_{k=0} {(-1)^k (1 \dot + q)^k x^k \over 
(1 -^{\hskip-.09truein \cdot} qx)^{k +1} }. 
\tag{16}
\end{equation*} 
By Euler's formula
\begin{equation*}
{1 \over (1 -^{\hskip-.09truein \cdot} x)^{k+1}} = \sum^\infty_{n=0}
 \bigg[{k+n \atop n}\bigg] x^n, \tag{17}
\end{equation*}
the RHS of (16) can be transformed into
\begin{equation*}
\sum^\infty_{k=0} (-x)^k (1 \dot + q)^k \sum^\infty_{n=0} \bigg[{k+n \atop k} \bigg] (qx)^n 
 = \sum^\infty_{N=0} x^N \sum^N_{k=0} \bigg[{N \atop k}\bigg] (-1)^k (1 \dot + q)^k q^{N-k} . 
\tag{18}
\end{equation*}
Thus, we need to verify that 
\begin{equation*}
(-1)^N = \sum^N_{k=0} \bigg[{N \atop k}\bigg] (-1)^k (1 \dot + q)^k q^{N-k}, \tag{19}
\end{equation*}
which is equivalent to
\begin{equation*}
1 = \sum^N_{k=0} \bigg[{N \atop k}\bigg] (-q)^k (1 \dot + q)^{N-k}, \tag{20}
\end{equation*}
which is true because [21, formula (2.10)]
\begin{equation*}
\sum^N_{k=0} \bigg[{N \atop k}\bigg] a^k (b \dot + v)^{N -k} = \sum^N_{k=0} \bigg[{N \atop 
k} \bigg] b^k (a \dot + v)^{N-k}  \tag*{\qed\ (21)}
\end{equation*}
\renewcommand{\qed}{}
\end{proof}
\begin{remark}{22} 
Formula (21) is a particular case $\{u = 0\}$ of the formula
\begin{equation*}
\sum^N_{k=0} \bigg[{N \atop k}\bigg] (a \dot + u)^k (b \dot + v)^{N-k} = \sum^N_{k=0} 
\bigg[{N\atop k}\bigg] (b \dot + u)^k (a \dot + v)^{N-k}. \tag{23}
\end{equation*}
\end{remark}
\begin{remark}{24}
What is the {\it{origin}} of formula (14)?  I don't know; I found this 
formula experimentally.
\end{remark}

Having quantised formula (12), we now reverse-engineer it.  Formula (11) has as a $q$-analog 
the easily verifyable equality
\begin{equation*}
{d \over d_qx} \bigg({x^k \over (1 -^{\hskip-.09truein \cdot} \alpha x)^k}\bigg) =
 [k] {x^{k-1} \over (1 -^{\hskip-.09truein \cdot} \alpha x)^{k+1}}, \tag{25}
\end{equation*}
where
\begin{equation*}
{d \over d_q x} (f (x)) = {f (qx) - f(x) \over qx - x} 
\end{equation*}
is the $q$-derivative.  Thus, formula (14) results from $q$-differentiating the equality
\begin{equation*}
\sum^\infty_{s=0}  {x^{2s+1} \over [2s+1]} = \sum^\infty_{k=0} {(-1)^k  (1 \dot + q)^k 
\over [k+1]} {x^{k+1} \over (1 -^{\hskip-.09truein \cdot} x)^{k+1}}, \tag{26}
\end{equation*}
both sides of which vanish when $x=0$.  This is our $q$-analog of formula (10).  Now, replace 
in formula (26) $x$ by $qx$ and then multiply both sides by $(1-x)^{-1}$, resulting in
\begin{equation*}
{1 \over 1-x} \sum^\infty_{s=0} {x^{2s+1} q^{2s+1} \over [2s+1]} = \sum^\infty_{k=0} 
{(-1)^{k-1} (1 \dot + q)^{k-1} q^k x^k \over (1 -^{\hskip-.09truein \cdot} x)^{k+1} [k]}. \tag{27}
\end{equation*}
By formula $(6c)$, the LHS of this identity equals to 
\begin{equation*}
\sum^\infty_{N=1} g_N x^N, \tag{28}
\end{equation*}
where
\begin{equation*}
g_N = \sum_{0 \leq 2s+1\leq N} {q^{2s+1} \over [2s+1]}; \tag{29}
\end{equation*}
the RHS of formula (27) can be transformed as in $(6e)$:
\begin{equation*}
\sum^\infty_{k=1} {(-1)^{k-1} (1 \dot + q)^{k-1} (qx)^k \over [k]} \sum^\infty_{n=0} 
\bigg[{k+n \atop n}\bigg] x^n = \sum^\infty_{N=1} x^N \sum^N_{k=1} {(-1)^{k-1} (1 \dot + q) 
^{k-1} \over [k]}  \bigg[{N \atop k}\bigg] q^k. \tag{30}
\end{equation*}
Thus,
\begin{equation*}
g_N = \sum_{0 \leq 2s+1\leq N} {q^{2s+1} \over [2s+1]} = \sum^N_{k=1} \bigg[{N \atop k}\bigg]
{(-1)^{k-1} (1 \dot + q)^{k-1} \over [k]} q^k. \tag{31}
\end{equation*}

This is the desired $q$-analog of formula (2), and we would have been entirely justified were we to 
stop right here, but it would be silly, not to mention quite unfair to Ramanujan who had 
uncovered for us the interesting formula (2).  So let's see what else lurks around our 
calculations.

Let's start with formula (14).  It can be rewritten as 
\begin{equation*}
{1 \over 1-x} + {1 \over 1+x} = \sum^\infty_{k=0} {(-x)^k (1 \dot + 1)^{k+1} \over 
(1 -^{\hskip-.09truein \cdot} x)^{k+2}}, \tag{32}
\end{equation*}
and this a particular case of the more general formula
\begin{equation*}
{1 \over 1-x} + {a \over 1 + ax} = \sum^\infty_{k=0} {(-x)^k  
(a \dot + 1)^{k+1} \over 
(1 -^{\hskip-.09truein \cdot} x)^{k+2}}, \tag{33}
\end{equation*}
which is equivalent to the formula
\begin{equation*}
1 + (-1)^N a^{N+1} = \sum^N_{k=0} \bigg[{N+1 \atop k+1}\bigg] (-1)^k (a \dot + 1)^{k+1}, 
\tag{34}
\end{equation*}
which is true by formula (23).

Formula (33), is turn, is the $q$-derivative of the identity
\begin{equation*}
-\Log (1-x) + \Log(1+ax) = \sum^\infty_{k=0} {(-1)^k x^{k+1} (a \dot + 1)^{k+1} \over 
[k+1] (1 -^{\hskip-.09truein \cdot} x)^{k+1}}, \tag{35}
\end{equation*}
which holds true in view of formula (25) and the vanishing of both sides of the equality 
(35) for $x=0$; here 
\begin{align*}
\Log (1+z) &= \Log (1 + z;q) = \sum^\infty_{k=0} {(-1)^k z^{k+1} \over [k+1]} \tag{36a}
\\
&= \int^z_0 {d_qt \over 1+t} \tag{36b}
\end{align*}
is the $q-$logarithm.  Replacing in formula (35) $a$ by $b$ and subtracting, we get
\begin{equation*}
\Log (1+ax) - \Log (1 + bx) = \sum^\infty_{k=0} {(-1)^k x^{k+1} ((a \dot + 1)^{k+1} 
 - (b \dot + 1)^{k+1} ) \over [k+1] (1 -^{\hskip-.09truein \cdot} x)^{k+1}}, \tag{37}
\end{equation*}
a rather peculiar formula.

Next, let's look again at Ramanujan's formula (2).  In his later years, Ramanujan would 
sometime work out on a  new formula without fully developing all the related consequences of 
those he had already found, but this was definitely not so in his younger years, and Chapter 
4 of his 2$^{nd}$ Notebook, where formula (1) is located, was recorded at the beginning of his 
mathematical journey.  Since this entry is not related to anything else in Chapter 4, we must 
conclude that either formula (2) is an isolated one and not a part of a general pattern; or 
that there {\it{exists}} such a pattern but Ramanujan had missed it; or that he didn't look 
for such a pattern at all, having more interesting problems to occupy himself with at the 
time.  We shall see in a moment that such a pattern {\it{does exist}}, and since it's 
inconceivable that Ramanujan could have missed it, the third alternative is almost certainly 
true.

Fix a positive integer $L \geq 2$.  Denote by
\begin{equation*}
G_{n|L} = \sum_{{{1 \leq k \leq n } \atop k \not\equiv 0 \ (\bmod \ L)}} {1 \over k} . \tag{38}
\end{equation*}
Then
\begin{align*}
G_L (x) &= \sum^\infty_{N=1} G_{N|L} x^N = \sum^\infty_{{{k=1 \atop k \not\equiv 0 \ (\bmod \ L)}}} 
{1 \over k} \sum^\infty_{s \geq 0} x^{k+s} 
\\
&={1 \over 1-x} \bigg(\sum^\infty_{k=1} {x^k \over k} - \sum^\infty_{k=1} {x^{kL} \over 
kL} \bigg) = {1 \over 1-x} \bigg(- \log (1 - x) + {1 \over L} \log (1 - x^L) \bigg) 
\\
&={1 \over L (1-x)} \log \bigg({1-x^L \over (1-x)^L} \bigg) = {1 \over L(1-x)} 
\log \bigg(1+ {1-x^L - (1-x)^L \over (1-x)^L} \bigg). \tag{39}
\end{align*}
The case $L=2$ is the Ramanujan case $(6a-d)$.  The next case is $L=3$:
\begin{align*}
G_3 (x) &= {1 \over 3(1-x)} \log \bigg({1 - x^3 \over (1-x)^3 } \bigg) = {1 \over 3(1-x)} 
\log \bigg({1+x+x^2 \over (1-x)^2} \bigg) 
\\
 &= {1 \over 3(1-x)} \log \bigg(1 + {3x \over (1-x)^2} \bigg) = \sum^\infty_{k=1} {3^{k-1} x^k 
\over (1-x)^{2k+1}} {(-1)^{k-1} \over k} 
\\
&= \sum^\infty_{k=1} {(-3)^{k-1} x^k \over k} \sum^\infty_{n=0} \bigg({2k+n \atop 2k} \bigg) 
x^n = \sum^\infty_{N=1} x^N \sum^N_{k=1} {(-3)^{k-1} \over k} \bigg({N+k \atop 2k} \bigg). 
\tag{40}
\end{align*}
Thus, 
\begin{equation*}
G_{3|N} = \sum_{{1 \leq k \leq N \atop k \not\equiv 0 \ (\bmod \ 3)}} {1 \over k} = \sum^N
_{k=1} {(-3)^{k-1}  \over k} \bigg({N + k \atop 2k} \bigg). \tag{41}
\end{equation*}

For $L > 3$, formulae become more complex, and we shan't pursue them further.

We shall return to formula (41) later on.  Let's now look at the following formula, 
similar $-$ in spirit if not in substance $-$ to the Ramanujan formula (2):
\begin{equation*}
\sum^N_{k=1} \bigg({N \atop k}\bigg) {(-1)^{k-1} \over k} - \sum^N_{k=1} {(-1)^{k-1} \over k} 
= \sum^{\lfloor N/2 \rfloor}_{k=1} {1 \over k}. \tag{42}
\end{equation*}
To prove this formula, discovered experimentally, we multiply both sides of it by $x^N$ 
and then sum on $N \in \N$.  We get:
\begin{align*}
\text{1) }\ & \sum^\infty_{N=1} x^N \sum^N_{k=1} \bigg({N \atop k}\bigg) {(-1)^{k-1} \over k} 
= \sum^\infty_{k=1} {(-1)^{k-1} \over k} x^k \sum^\infty_{N=k} \bigg({N \atop k}\bigg) 
x^{N-k} 
\\
&\qquad = \sum^\infty_{k=1} {(-1)^{k-1} x^k \over k} {1 \over (1-x)^{k+1}} 
= {1 \over 1-x} \sum^\infty_{k=1} {(-1)^{k-1} \over k} \bigg({x \over 1-x}\bigg)^k 
\\
&\qquad = {1 \over 1-x} \log \bigg({1+ {x \over 1-x}} \bigg) 
= {1 \over 1-x} \log \bigg({1 \over 1-x}\bigg) = - 
{1 \over 1-x} \log (1-x);  \tag{43a}
\\[1ex]
\text{2) }\ & \sum^\infty_{N=1} x^N \sum^N_{k=1} {(-1)^{k-1} \over k} = \sum^\infty_{k=1} 
{(-1)^{k-1} \over k} \sum^\infty_{N=k} x^N 
\\
&\qquad = \sum^\infty_{k=1} {(-1)^{k-1} \over k} 
{x^k \over 1-x} = {1 \over 1-x} \log (1+x); \tag{43b}
\\[1ex]
\text{3) }\ &\sum^\infty_{N=1} x^N \sum^{\lfloor N/2 \rfloor}_{k=1} {1 \over k} = \sum^\infty_{k=1} 
{1 \over k} \bigg(\sum^\infty_{s=0} x^{2k+2s} + x^{2k+1+2s}\bigg) = \sum^\infty_{k=1} 
{x^{2k} \over k} (1 + x) {1 \over 1-x^2} 
\\
&\qquad = {-1 \over 1-x} \log (1-x^2). \tag{43c}
\end{align*}
Thus, formula (42) is equivalent  to the equality 
\begin{equation*}
\log (1-x) + \log (1+x) = \log (1-x^2), \tag{44}
\end{equation*}
which is obviously true.

Formula (42) consists of {\it{three}} different sums, so its quantization is quite likely 
to be hugely nonunique, and probably the less evident the more interesting.  The reader may 
wish to try to quantize this formula first before reading any further. 
\begin{proposition}{45}
We have:
\begin{equation*}
\sum^N_{k=1} \bigg[{N \atop k}\bigg] {(-1)^{k-1} q^{({k+1 \atop 2})} \over [k]} = 
\sum^N_{k=1} {(-1)^{k-1} q^k \over [k]} + \sum^{\lfloor N/2 \rfloor}_{k=1} {q^{2k} \over 
[2k]/2}. \tag{46}
\end{equation*}
\end{proposition}
\begin{proof}
Multiplying both sides of formula (46) by $x^N$ and summing on $N \in \N$, 
we find:
\begin{align*}
\text{1) }\ &\sum^\infty_{N=1} x^N \sum^N_{k=1} \bigg[{N \atop k}\bigg] {(-1)^{k-1} q^{{({k+1 \atop 
2})}} \over [k]} = \sum^\infty_{k=1} {(-1)^{k-1} q^{{({k+1 \atop 2})}} \over [k]} 
\sum^\infty_{N=k} \bigg[{N \atop k}\bigg] x^N 
\\
&\qquad = \sum^\infty_{k=1} {(-1)^{k-1} q^{{({k+1 \atop 2})}} \over [k]} {x^k \over 
(1 -^{\hskip-.09truein \cdot} x)^{k+1}} = {1 \over 1-x} \sum^\infty_{k=1} 
{(-1)^{k-1} q^{({k \atop 2})} (qx)^k \over [k] (1 -^{\hskip-.09truein \cdot} qx)^k}; 
\tag{47a}
\\
\text{2) }\ &\sum^\infty_{N=1} x^N \sum^N_{k=1} {(-1)^{k-1} q^k \over [k]} = {1 \over 1-x} 
\sum^\infty_{k=1} {(-1)^{k-1} (qx)^k \over [k]}; \tag{47b} 
\\
\text{3) }\ &\sum^\infty_{N=1} x^N \sum^{\rfloor N/2 \lfloor}_{k=1} {q^{2k} \over [2k]/2} = 
{1 \over 1-x} \sum^\infty_{k=1} {(qx)^{2k} \over [2k]/2}, \tag{47c}
\end{align*}
where we used the two handy formulae:
\begin{gather*}
\sum^\infty_{N=1} x^N \sum^{N}_{k=1} \varphi_k = {1 \over 1-x} \sum^\infty_{k=1} \varphi_k 
x^k, \tag{48a}
\\
\sum^\infty_{N=1} x^N \sum^{\lfloor N/2 \rfloor}_{k=1} \psi_k = {1 \over 1-x} 
\sum^\infty_{k=1} \psi_k x^{2k}. \tag{48b}
\end{gather*}
Formulae (47) convert the identity (46) into
\begin{equation*}
\sum^\infty_{k=1} {(-1)^{k-1} t^k q^{({k \atop 2})} \over [k] 
(1 -^{\hskip-.09truein \cdot} t)^k } = \sum^\infty_{k=1} {(-1)^{k-1} t^k \over [k]} + 
\sum^\infty_{k=1} {t^{2k} \over [2k]/2}, \tag{49}
\end{equation*}
where $t=qx$.  Each side of the equality (49) vanishes for $t=0$, so it's equivalent   
to the $q$-derivative ${{\dps{d \over d_q t}}} $ of it:

\begin{equation*}
\sum^\infty_{k=1} {(-t)^{k-1} q^{({k \atop 2})} \over (1 -^{\hskip-.09truein \cdot}t)^{k+1}} 
= \sum^\infty_{k=1} (-t)^{k-1} + 2 \sum^\infty_{k=1} t^{2k-1}, \tag{50}
\end{equation*}
where we used formula (25).  Since the RHS of formula (50) is
\begin{equation*}
\sum_{k \ {\rm odd} \ > 0}  (-t)^{k-1} + \sum_{k \ {\rm odd} \ > 0} t^k = \sum^\infty_{k=0} t^k = 
{1 \over 1-t}, 
\end{equation*}
we need to verify that
\begin{equation*}
\sum^\infty_{k=0} {(-t)^k q^{({k+1 \atop 2})} \over (1 -^{\hskip-.09truein \cdot} t)^{k+2} } 
= {1 \over 1-t}, \tag{51}
\end{equation*} 
and this is formula (33) for $a=0$.
\end{proof}
\begin{remark}{52}
The RHS of formula (49) is easily seen to be ${\dps{\sum^\infty
_{k=1} {t^k \over [k]}}}$.  This suggests that similar simplifications are possible also in 
other formulae figuring in the Proof of the identify (46); and indeed, the RHS of formula 
(46) becomes, on closer inspection, just ${\dps{\sum^N_{k=1} \ {q^k \over [k]}}}$.  Identity 
(46) then turns into 
\begin{equation*}
\sum^N_{k=1} \bigg[{N \atop k}\bigg] {(-1)^{k-1} q^{{({k+1 \atop 2})}} \over [k]} = 
\sum^N_{k=1} {q^k \over [k]}. \tag{53}
\end{equation*}
Formula (53) had appeared as Problem $\#$ 6407 in the  {\it{American Mathematical Monthly}}.
\end{remark}

The reader may wonder whether there is some common theme uniting under the single roof all the 
identities we have met so far.  At least in the classical case $q=1$, I think the common theme 
is the generating functions which are linear in logarithms, with coefficients that are 
rational functions themselves.

The reader may also wonder what the deal is with all these identities anyway; 
{\it{why bother}}?  Apart from showing the proper respect to Ramanujan, that is.  Perhaps Littlewood's 
candid observation [14, \ p. 103] is as good an explanation as any:  ``Mathematics is a  dangerous 
profession; an appreciable proportion of us go mad..."

We now return to formula (41).  In deriving it, we have gone through the calculation
\begin{gather*}
\sum_{k \not\equiv 0 \ (\bmod \ 3)} {x^k \over k} = \sum^\infty_{k=1} {x^k \over k} - 
\sum^\infty_{k=1} {x^{3k} \over 3k} = - \log (1-x) + {1 \over 3} \log (1-x^3) 
\\
\qquad = {1 \over 3} \log \bigg({1-x^3 \over (1-x)^3}\bigg) = {1 \over 3} \log \bigg( 1+ {3x \over 
(1-x)^2} \bigg) = \sum^\infty_{k=1} {(-1)^{k-1} 3^{k-1} x^k \over k (1 - x)^{2k}}, 
\end{gather*}
so that
\begin{equation*}
\sum^\infty_{k=1} {x^k \over k} - \sum^\infty_{k=1} {x^{3k} \over 3k} = \sum^\infty_{k=1} 
{(-1)^{k-1} 3^{k-1} x^k \over k (1-x)^{2k}}. \tag{54}
\end{equation*}
It is this formula we shall now quantize.
\begin{proposition}{55} 
\begin{equation*}
\sum_{k \not\equiv 0 \ (\bmod \ 3)} {x^k \over [k]} = \sum^\infty_{k=1} 
{(-1)^{k-1} <3^{k-1}> q^{{-({k \atop 2})}} x^k \over [k] (1 -^{\hskip-.09truein \cdot} q^{-k} 
x)^{2k}}, \tag{56}
\end{equation*}
where
\begin{equation*}
<3^n> = <3^n>_{\dps{q}} = \Pi^n_{k=1} [3]_{{\dps{q}}^{k}} , 
\ \ \ n \ \in \N; \ \ \ <3^0> = 1. \tag{57}
\end{equation*}
\end{proposition}
\begin{proof}
Each side of the identity (56) vanishes for $x=0$, so we apply 
${\dps{d \over d_qx}}$ to it to make it simpler-looking.  We get:
\begin{gather*}
{d \over d_q x} (LHS) = \sum^\infty_{k=1} x^{k-1} - \sum^\infty_{k=1} x^{3k-1} = 
{1 \over 1-x} - {x^2 \over 1-x^2} = {1+x \over 1-x^3}, \tag{58a}
\\
{d \over d_q x} (RHS) = \sum^\infty_{k=1} {(-x)^{k-1} <3^{k-1}> q^{-{({k \atop 2})}} 
(1 +x) \over (1 -^{\hskip-.09truein \cdot} q^{-k} x)^{2k+1}}, \tag{58b}
\end{gather*}
because, as is easily verified, 
\begin{equation*}
{d \over d_q x} \bigg({x^k \over (1 -^{\hskip-.09truein \cdot} \alpha x)^{2k}} \bigg) 
= [k] {x^{k-1} (1 + \alpha q^k x) \over (1 -^{\hskip-.09truein \cdot} \alpha x)^{2k+1} }. 
\tag{59}
\end{equation*}
We thus arrive at the identify
\begin{equation*}
{1 \over 1-x^3} = \sum^\infty_{k=0} {(-x)^k <3^k>q^{- {({k+1 \atop 2})}} \over 
(1 -^{\hskip-.09truein \cdot} q^{-k-1} x) ^{2k+3}}. \tag{60}
\end{equation*}
Denote
\begin{equation*}
S_N = \sum^N_{k=0} {(-x)^k <3^k>q^{-{({k+1 \atop 2})}} \over (1 -^{\hskip-.09truein \cdot}
q^{-k-1}x) ^{2k+3} } . \tag{61}
\end{equation*}
We shall prove that
\begin{equation*}
(1+x+x^2) S_N = {1 \over 1-x} + {(-1)^N x^{N+1} <3^{N+1}> q^{-({N+2 \atop 2})} \over 
(1 -^{\hskip-.09truein \cdot} q^{-N-1} x)^{2N+3}} . \tag{62}
\end{equation*}
Formula (60) then results from formula (62) when $N \rightarrow \infty$. 

To prove formula (62) we use induction on $N$.  For $N=0$, formula (60) is obviously true.  
Now, denote the RHS of formula (62) by $R_N$.  Then
\begin{align*}
R_N - R_{N-1} &= {(-1)^N x^{N+1} <3^{N+1}> q^{-({N+2 \atop 2})} \over 
(1-^{\hskip-.09truein \cdot}
q^{-N-1} x)^{2N+3}} - {(-1)^{N-1} x^N <3^N>q^{-({N +1 \atop 2})} \over 
(1 -^{\hskip-.09truein \cdot}
q^{-N} x)^{2N+1} }  
\\
&= {(-x)^N <3^N>q^{-({N+1 \atop 2})} \over (1-^{\hskip-.09truein \cdot} q^{-N-1} x)^{2N+3}} 
\Delta_N = (S_N - S_{N-1}) \Delta_N, \tag{63}
\end{align*}
where
\begin{align*}
\Delta_N &=  xq^{-N-1} [3]_{q^{N+1}} + (1 - q^{-N-1} x) (1 - q^{N+1} x) = 
\\
&= x (q^{-N-1} + 1 + q^{N+1}) + (1 - (q^{-N-1} + q^{N+1}) x + x^2) = 1 +x+x^2. \tag{64}
\end{align*}
Thus,
\begin{equation*}
R_N - R_{N-1} = (1 + x + x^2) (S_N - S_{N-1}). \tag*{\qed}
\end{equation*}
\renewcommand{\qed}{}
\end{proof}
\begin{remark}{65}
Formula (60) can be rewritten as
\begin{equation*}
\sum^N_{k=0} (-1)^k q^{-({k+1 \atop 2})} <3^k> q^{-(k+1)(N-k)} \bigg[{N+k+2 \atop 2k+2}
\bigg] = \cases{1, & $N \equiv 0 \ (\bmod \ 3)$, \cr
0, & $N \not\equiv 0 \ (\bmod \ 3)$. \cr} \tag{66}
\end{equation*}
This formula is invariant w.r.t.\ the change of $q$ into $q^{-1}$. 
\end{remark}

>From our findings so far, let's collect together formulae (33)$|_{a=0},$ (14), and (60):
\begin{gather*}
{1 \over 1-x} = \sum^\infty_{k=0} {(-x)^k \over (1 - x)^{k+2}} = \sum^\infty_{k=0} 
{(-x)^k q^{({k+1 \atop 2})} <1^k> \over (1 -^{\hskip-.09truein \cdot} x)^{k+2}}, \tag{67.1}
\\
{1 \over 1-x^2} = \sum^\infty_{k=0} {(-x)^k 2^k \over (1 - x)^{k+2}} = \sum^\infty_{k=0} 
{(-x)^k <2^k> \over (1 -^{\hskip-.09truein \cdot} x)^{k+2}}, \tag{67.2}
\\
{1 \over 1-x^3} = \sum^\infty_{k=0} {(-x)^k 3^k \over (1 - x)^{2k+3} } = \sum^\infty_{k=0} 
{(-x)^k q^{-({k+1 \atop 2})} <3^k> \over (1 -^{\hskip-.09truein \cdot} q^{-k-1} x)^{2k+3} }. 
\tag{67.3}
\end{gather*}

We see that we have gotten ourselves the highest mathematical prize possible:  a fruitful 
definition, for formulae (67) suggest that, in {\it{some circumstances}}, a proper quantum 
version of the number $x^n$ is not the straightforward $[x^n]_q$ but a slightly devious 
\begin{equation*}
<x^n> = <x^n>_{{\dps{q}}} = \Pi^n_{k=1} 
[x]_{{\dps{q^{k}}}} , \ \ \ n \in \N; \ \ <x^0 > = 1. \tag{68}
\end{equation*}
It is not immediately clear exactly how useful this definition really is; we shall examine 
some supporting evidence in a little while.  At the moment, let's try to generalize the 
definition (68) for the case when $n$ is arbitrary (real number, complex number, formal 
parameter, ...) rather than a nonnegative integer.  (So that, e.g., a related $q$-dzeta function 
could be defined.  E.g., $\sum^\infty_{n=1} {q^{\alpha (n)} \over <n^s>}$, where $\alpha (n)$ is 
a quadratic polynomial in $n$. See [15,16].)  Since
\begin{equation*}
<x^n> = \Pi^n_{k=1} [x]_{{\dps{q}}^{k}} = \Pi^n_{k=1} {1 - q^{xk} \over 1-q^k} = 
{(q^x; q^x)_n \over (q; q)_n}, \tag{69}
\end{equation*}
where
\begin{equation*}
(a; b)_n = \Pi^{n-1}_{s=0} (1 - ab^s), \ \ \ n \in \N; \ \ \ (a; b)_0 = 1, \tag{70}
\end{equation*}
and
\begin{equation*}
(a; b)_n = {(a; b)_\infty \over (ab^n;b)_\infty}, \tag{71}
\end{equation*}
the latter formula {\it{defining}} the symbol $(a;b)_n$ for arbitrary $n$ (with $|q|<1$ 
being assumed for convergence), we have:
\begin{equation*}
<x^n> = {(q^x; q^x)_\infty \over (q^{(n+1)x}; q^x)_\infty} {(q^{n+1}; q)_\infty \over (q; q)
_\infty}. \tag{72}
\end{equation*}
This is our desired definition.  We can also define 
\begin{equation*}
<x^\infty > = \Pi^\infty_{k=1} [x]_{{\dps{q}}^k} = {(q^x; q^x)_\infty \over (q; q)_\infty}, 
\tag{73}
\end{equation*}
so that
\begin{equation*}
<x^n> = <x^\infty> {(q^{n+1}; q)_\infty \over (q^{(n+1)x}; q^x)_\infty}. \tag{74}
\end{equation*}
Also, 
\begin{equation*}
<x^n> = [x]^n {\Gamma _{q^{x}}{(n+1}) \over \Gamma_{q}{(n+1)} }, \tag{75}
\end{equation*}
where
\begin{equation*}
\Gamma _{\dps{q}} (n+1) = (1 - q)^{-n} {(q; q)_\infty \over (q^{n+1}; q)_\infty} \tag{76a}
\end{equation*}
is the $q$-Gamma function; for positive integers $n$,
\begin{align*}
\Gamma_{\dps{q}} (n+1) &= [n]_q! = \Pi^n_{k=1} [k]_{\dps{q}} = \Pi^n_{k=1} {1-q^k \over 1-q} 
\\
&= (1 - q)^{-n} (q; q)_n = (1-q)^{-n} {(q;q)_\infty \over (q^{n+1}; q)_\infty}. \tag{76b}
\end{align*}
 
The objects $<x^\infty >$ are interesting in their own right, although we won't pursue them 
here.  As just one example, if $L \geq 2$ is an integer, then
\begin{equation*}
{<(2L)^\infty > \over < L^\infty > < 2^\infty> } = \Pi_{n \not\equiv 0 \ (\bmod \ L)} 
(1 + q^n)^{-1}. \tag{77}
\end{equation*}

We now look at a few classical situations where the numbers $2^n$ and $x^n$ enter into the 
picture, and examine whether the quantization recipe $x^n \mapsto \ <x^n>$ (68) provides agreeable 
results or not.

Let us start with the simple geometric progression:
\begin{equation*}
{1 \over 2} + {1 \over 4} + {1 \over 8} + ... = \sum^\infty_{k=0} {1 \over 2^{k+1}} = 1. 
\tag{78}
\end{equation*}
If we take formula (67.1), multiply both parts of it by $1-x$, and then set $x=-1$, we get
\begin{equation*}
1 = \sum^\infty_{k=0} {q^{({k+1 \atop 2})} \over <2^{k+1}>}. \tag{79}
\end{equation*}
The substitution $x=-1$ is of course contingent upon the series (67.1) being convergent 
for $x=-1$.  Alternatively, induction on $N$ shows that
\begin{equation*}
\sum^N_{k=0} {q^{({k+1 \atop 2})} \over <2^{k+1}>} = 1 - {q^{({N+2 \atop 2})} \over 
<2^{N+1}>}. \tag{80}
\end{equation*}
 
Our next subject is the Euler transformation of series [20, p. 244]:
\begin{equation*}
\sum^\infty_{k=0} (-1)^k a_k = \sum^\infty_{\ell=0} {(-1)^\ell (\Delta^\ell a)_0 \over 
2^{\ell +1}}, \tag{81}
\end{equation*}
where $\Delta$ is the difference operator:
\begin{equation*}
(\Delta^0 a)_k = a_k; \ \ (\Delta^{\ell +1} a)_k = (\Delta^\ell a)_{k+1} - (\Delta^\ell a)_k, 
\ \ \ell \in \Z_+. \tag{82}
\end{equation*}
We quantize these formulae thusly:
\begin{gather*}
(\Delta^0 a)_k = a_k; \ \  (\Delta^{\ell+1} a)_k = (\Delta^\ell a)_{k+1} - q^\ell (\Delta^\ell 
a)_k, \ \ \ell \in \Z_+, \tag{83}
\\[1ex]
\sum^\infty_{k=0} (-q)^k a_k = \sum^\infty_{\ell=0} {(-q)^\ell (\Delta^\ell a)_0 \over 
<2^{\ell +1} >}. \tag{84}
\end{gather*}
To prove  the identity (84) we start off with the easily proven by induction on $\ell$ 
formula
\begin{equation*}
(\Delta^\ell a)_k = \sum^\ell_{s=0} a_{k+\ell-s} (-1)^s q^{({s \atop 2})} \bigg[{\ell \atop s}
\bigg]. \tag{85}
\end{equation*} 
In particular,
\begin{equation*}
(\Delta^\ell a)_0 = \sum^\ell_{s=0} a_{\ell-s} (-1)^s q^{({s \atop 2})} \bigg[{\ell \atop s}
\bigg]. \tag{86}
\end{equation*}
Therefore, for the RHS of formula (84) we get:
\begin{gather*}
\sum^\infty_{\ell =0} (-q)^\ell {(\Delta^\ell a)_0 \over <2^{\ell +1}>} = \sum_{\ell, s} 
a_{\ell -s} (-1)^s q^{({s \atop 2})} \bigg[{\ell \atop s}\bigg] (-q)^\ell {1 \over <2^{\ell +1} 
> } 
\\
\qquad = \sum^\infty_{k=0} a_k (-q)^k \sum^\infty_{s=0} q^{({s \atop 2})} \bigg[{k+s \atop s}\bigg]
q^s {1 \over (1 \dot + q)^{k+s+1}} \ {\rm{[by}} \ (88)] = \sum^\infty_{
k=0} a_k (-q)^k. 
\end{gather*}
\begin{proposition}{87} 
\begin{equation*}
\sum^\infty_{s=0} q^{({s \atop 2})} \bigg[{k+s \atop s}\bigg] q^s 
 {1 \over (1 \dot +
q)^{k+s+1}} = 1, \ \ \ \forall k \in \Z_+. \tag{88}
\end{equation*}
\end{proposition}
\begin{proof}
By formula (17), the LHS of formula (88) is
\begin{gather*}
\sum_{s,n} q^{({s \atop 2})} \bigg[{k+s \atop s}\bigg] q^s \bigg[{k+s+n \atop n}\bigg] 
(-q)^n = \sum_{s,n} \bigg[{k+s+n \atop k}\bigg] \bigg[{s+n \atop s}\bigg] q^{({s \atop 2})} 
q^s (-q)^n 
\\
\qquad = \sum_N \bigg[{N+k \atop k} \bigg] \sum_{s+n=N} \bigg[{s+n \atop s}\bigg] q^{({s \atop 2})} 
q^s (-q)^n = \sum_N \bigg[{N+k \atop k} \bigg] (-q \dot + q)^N 
\\
\qquad = \sum_N \bigg[{N+k \atop k} \bigg] 
\delta^N_0 = 1 \tag*{\qed}
\end{gather*}
\renewcommand{\qed}{}
\end{proof}

Knopp gives two examples of the Euler transformation, on p. 246, ibid.  In the first,
\begin{gather*}
a_k = {1 \over k+1}, \ \ (\Delta^\ell a)_k = {(-1)^\ell \ell! \over (k+1)... (k+\ell+1)} , \ \  
(\Delta^\ell a)_0 = {(-1)^\ell \over \ell +1}, \tag{89a}
\\[1ex]
\log 2 = \sum^\infty_{k=0} {(-1)^k \over k+1} = \sum^\infty_{\ell =0} {1 \over (\ell +1)2^{\ell+1}}, 
\tag{89b}
\end{gather*}
and in the second, 
\begin{gather*}
a_k = {1 \over 2k+1} , \ \ (\Delta^\ell a)_k = {(-1)^\ell (2 \ell)!! \over \Pi_{s=0}^\ell 
(2k+1+2s)} , \ \ (\Delta^\ell a)_0 = {(-1)^\ell (2 \ell)!! \over (2 \ell +1)!!}, \tag{90a}
\\[1ex]
{\pi \over 4} = \sum^\infty_{k=0} {(-1)^k \over 2k+1} = {1 \over 2} (1 + {1 ! \over 1 \cdot 3} 
+ {2! \over 1\cdot 3 \cdot 5} + {3! \over 1 \cdot 3 \cdot 5 \cdot 7} + ...) = \sum^\infty_{\ell =0} 
{(2\ell)!! \over 2^{\ell +1} (2 \ell +1)!!} . \tag{90b}
\end{gather*}
Quantum versions of these two examples are:
\begin{gather*}
a_k = {1 \over [k+1]}, \ \ (\Delta^\ell a)_k = {(-1)^\ell q^{k \ell + ({\ell + 1 \atop 2})} 
[\ell]! \over [k+1]... [k+\ell +1]} , \ \ (\Delta^\ell a)_0 = {(-1)^\ell q^{({\ell+1 \atop 2})} 
\over [\ell + 1]}, \tag{91a}
\\[1ex]
q^{-1} \Log ([2]; q)= \sum^\infty_{k=0} {(-q)^k \over [k+1]} = \sum^\infty_{\ell=0} 
{q^{({\ell +1 \atop 2})} \over [\ell +1] < 2^{\ell+1}>}, \tag{91b}
\end{gather*}
and
\begin{gather*}
a_k = {1 \over [2k+1]_{q^{1/2}} }, \ \ (\Delta^\ell a)_k = {(-1)^\ell q^{(k+1/2)\ell+ 
{({\ell \atop 
2})} } [2 \ell ]_{q^{1/2}} !! \over \Pi^\ell_{s=0} [2k+1+2s]_{q^{1/2}}}, 
\\
(\Delta^\ell a)_0 = {(-1)^\ell q^{{\ell \over 2}
+ {({\ell \atop 2})}} [2 \ell]_{q^{1/2}} !!
\over [2 \ell +1]_{q^{1/2}}!!}, \tag{92a}
\\
\bigg({\pi \over 4}\bigg)_{q^{1/2};q} =
\int^1_0 {d_{q^{1/2}}t \over {1+qt^2}}\\
\qquad\qquad\ \ {= \sum^\infty_{k=0}  
{(-q)^k \over [2k+1]_{q^{1/2}}}
= \sum^\infty_{\ell =0} {q^{{\ell \over 2}
+ {({\ell \atop 2})}} [2\ell]_{q^{1/2}}!!
\over <2^{\ell+1} >_q [2 \ell +1]_{q^{1/2}}!!}}. \tag{92b}
\end{gather*}


We next venture into the virgin-so-far territory of Quantum Number Theory (QNT).  (``What 
is QNT?''  Turn your imagination on and find out.)  Let $p$ be a 
prime and {\it{a}} be a positive integer coprime to $p$.  The (Small) Fermat Theorem (SFT) states
that 
\begin{equation*}
a^{p-1} \ \equiv 1 \ \ (\bmod \ p), \ \ \ (a, p)= 1, \tag{93}
\end{equation*}
and therefore
\begin{equation*}
a^p \equiv a \ \ (\bmod \ p), \tag{94}
\end{equation*}
with the latter equality being true for all $a$'s, not just those coprime to $p$.
\begin{proposition}{95}
If $a$ is coprime to a prime $p$ then
\begin{equation*}
<a^{p-1}> \equiv 1 \ \ (\bmod \ [p]), \ \ (a, p) = 1, \tag{96}
\end{equation*}
and 
\begin{equation*}
<a^p> \equiv a \ (\bmod \ [p]), \ \ (a, p) = 1. \tag{97}
\end{equation*}
If $a$ is not coprime to $p$, then
\begin{equation*}
<a^p> \equiv [a] \ (\bmod \ [p]), \ \ (a, p ) \not= 1. \tag{98}
\end{equation*}
(All congruences of QNT live in $\Z[q], \ \ \Z[q^{-1}, q], $ etc.)
\end{proposition}
\begin{proof}
Since
\begin{equation*}
(q-1) [p] = q^p-1 \equiv 0 \ (\bmod \ [p]), \tag{99}
\end{equation*}
we have 
\begin{equation*}
q^p \equiv 1 \ (\bmod \ [p]). \tag{100}
\end{equation*}
Now, 
\begin{equation*}
<a^{p-1}> = \Pi^{p-1}_{k=1} (1-q^{ak}) / \Pi^{p-1}_{k=1} (1-q^k). \tag{101}
\end{equation*}
Since, ($a, p)=1$, the exponents $\{a, 2a, ... , (p-1)a\}$ are, modulo $p$, just a 
permutation of the exponents $\{1, 2,...,p-1\}$.  
Therefore, the expressions $\{q^a, q^{2a}, ... , q^{(p-1)a} \}$ are, since $q^p \equiv 1 \ \ 
(\bmod \ [p])$, just a permutation of the expressions $\{q, q^2, ... , q^{p-1})$ modulo 
$[p]$.  This proves (96). 

Now, for {\it{any}} $a$, coprime to $p$ or not, 
\begin{equation*}
[a]_{{\dps{q}}^{p}} = 1+q^p + ... + q^{(a-1)p} \equiv a  \ (\bmod \ [p]), \tag{102}
\end{equation*}
so that, by (96), for {\it{a}} coprime to $p$, 
\begin{equation*}
<a^p> = <a^{p-1}> [a]_{{\dps{q}}^{p}} \equiv 1 \cdot a = a \ (\bmod \ [p]). \tag{103}
\end{equation*}
This proves (97).

Finally, if $a$ is divisible by $p$, 
\begin{equation*}
 a = p \ell, 
\end{equation*}
then 
\begin{equation*}
[a] = [p \ell] = [p] [\ell]_{{\dps{q}^{p}} }\equiv 0 \ (\bmod \ [p]), \tag{104}
\end{equation*}
so that
\begin{equation*}
<a^p> = \Pi^p_{k=1} [a]_{{\dps{q}}^{k}} \equiv 0 \ (\bmod  \ [p]). \tag{105}
\end{equation*}
This proves  (98).
\end{proof}
\begin{remark}{106}
Some quantum congruences have no classical analogs.  For example: 
\begin{proposition}{107}
For a prime $p$,
\begin{equation*}
(1 -^{\hskip-.09truein \cdot} q)^{p-1} \equiv p \ (\bmod \ [p]). \tag{108}
\end{equation*}
\end{proposition}
\begin{proof}
First, since $p$ is a prime, 
\begin{equation*}
\bigg[{p \atop k}\bigg] \equiv \left\{\matrix{0, & 0 < k < p & \cr
 \ \ \ & & (\bmod \ [p]). \cr
1, & k=0 \ {\rm{or}} \  p & \cr}\right. \tag{109}
\end{equation*}
Next, it's easy to verify that
\begin{equation*}
\bigg[{n \atop k}\bigg] = \sum^k_{s=0} \bigg[{n+1 \atop k-s}\bigg] (-1)^s q^{s(n-k)} q^{({s+1 
\atop 2})}. \tag{110}
\end{equation*}
Taking $n=p-1$ and using formula (109), we get
\begin{equation*}
\bigg[{p-1 \atop k}\bigg] \equiv (-1)^k q^{-{({k+1 \atop 2})}} \  (\bmod \ [p]). \tag{111}
\end{equation*}
Therefore, 
\begin{equation*}
(1 -^{\hskip-.09truein \cdot} q)^{p-1} = \sum^{p-1}_{k=0} (-q)^k \bigg[{p-1 \atop k}\bigg] 
q^{({k \atop 2})} \equiv \sum^{p-1}_{k=0} 1 = p \ (\bmod \ [p]) \tag*{\qed (112)}
\end{equation*}
\renewcommand{\qed}{}
\end{proof}

In exactly the same way, we obtain 
\begin{equation*}
(1-^{\hskip-.09truein \cdot} qx)^{p-1} \equiv {1 - x^p \over 1-x} \ (\bmod \ [p]), \tag{113}
\end{equation*}
and  while we are at it, let's notice that
\begin{equation*}
(1 \dot + x)^p \equiv 1+x^p \ (\bmod \ [p]), \tag{114}
\end{equation*}
a quantum version of the classical
\begin{equation*}
(1 + x)^p \equiv 1+x^p \ (\bmod \ p). \tag{115}
\end{equation*}
\end{remark}

Let now $m \in \N$ be arbitrary and not necessarily a prime.  Let $\varphi (m)$ be the 
total number of positive integers less than $m$ and coprime to $m$.  Euler's form of the 
Small Fermat Theorem states that for any {\it{a}} coprime to $m$, 
\begin{equation*}
a^{\varphi(m)} \equiv 1 \ (\bmod \ m), \ \ (a, m)=1. \tag{116}
\end{equation*}
\begin{proposition}{117}
Denote by $r_1, ... , r_{\varphi (m)}$ the complete set of 
positive integers less than $m$ that are coprime to $m$.  Then 
\begin{equation*}
\Pi^{\varphi(m)}_{i=1} [a]_{{\dps{q}}^{r_{i}}} \equiv 1 \ (\bmod \ [m]), \ \ \ (a, m)=1. \tag{118}
\end{equation*}
\end{proposition}
\begin{proof}
The LHS of the congruence (118) is 
\begin{equation*}
\Pi^{\varphi(m)}_{i=1} (1-q^{ar_{i}}) / \Pi^{\varphi(m)}_{i=1} 
(1-q^{r_{i}}). \tag{119}
\end{equation*}
Since $a$ is coprime to $m$, the set $\{ar_1, ..., ar_{\varphi(m)}\}$ is, modulo $m$, a 
permutation of the set $\{r_1, ... , r_{\varphi (m)}\}$.
\end{proof}

When $m=p$ is prime, formula (118) becomes formula (96).
\begin{corollary}{120}
For any $r=r_j$, the congruence 
\begin{equation*}
[a]_{q^{r}} x \equiv b \ (\bmod \ [m]) \tag{121}
\end{equation*}
has a solution
\begin{equation*}
x = b \Pi_{i \not= j} [a]_{q^{r_{i}}} . \tag{122}
\end{equation*}
\end{corollary}

We next turn to the classical Wilson theorem:  if $p$ is prime then
\begin{equation*}
(p-1)! \equiv -1 \ (\bmod \ p). \tag{123}
\end{equation*}
\begin{proposition}{124}
Let $p$ be a prime.  For each integer $a, \ 
 0 < a < p$, let 
$\bar a$ be the unique solution of the congruence 
\begin{equation*}
a \bar a \equiv 1 \ (\bmod \ p), \ \ \ 0 <  \bar a < p. \tag{125}
\end{equation*}
(Thus, $\bar a = a^{-1} \ {\rm{in}} \ \Z_p.)$  Then
\begin{align*}
\Pi^{p-1}_{a=1} [a]_{{\dps{q}}^{\bar a -1}} &\equiv - q^{-1} \ (\bmod \ [p]) \tag{126a}
\\
&\equiv [p-1] \ (\bmod \ [p]). \tag{126b}
\end{align*}
\end{proposition}
\begin{proof}
The congruence
\begin{equation*}
x^2 \equiv 1 \ (\bmod \ p) \tag{127}
\end{equation*}
has exactly two solutions:
\begin{equation*}
x=1, \ \ \ x=p-1, \tag{128}
\end{equation*}
so that
\begin{equation*}
\bar 1 = 1, \ \ \ \overline{p-1} = p-1. \tag{129}
\end{equation*}
Otherwise, 
\begin{equation*}
a \not= \bar a, \ \ \ 1 < a, \bar a < p-1. \tag{130}
\end{equation*}
For each such pair $a \not= \ \bar a, $ we have:
\begin{align*}
[a]_{{\dps{q}}^{\bar a-1}} [\bar a]_{{\dps{q}}^{a-1}} &= {1 - q^{(\bar a-1)a} 
\over 1-q^{\bar a-1}} 
{1-q^{(a-1)\bar a} \over 1-q^{a-1}} 
\\
&= {1 - q^{a \bar a-a} \over 1-q^{\bar a-1}} \cdot {1 - q^{a \bar a- \bar a}
 \over 1-q^{a-1}} \ \  
[{\rm{by}} \ (125), (100)] \equiv {(1-q^{1-a}) (1-q^{1-\bar a}) \over 
(1-q^{\bar a-1}) (1-q^{a-1})} 
\\
&= q^{1-a} q^{1-\bar a} = q^{2- (a+\bar a)} \ (\bmod \ [p]). \tag{131}
\end{align*}
Therefore, for all such ${{\dps{p-3 \over 2}}}$ pairs combined,
\begin{equation*}
\Pi^{p-2}_{a=2} [a]_{{\dps{q}}^{\bar a-1}} \equiv q^\sigma \ \ (\bmod \ [p]), 
\end{equation*}
where
\begin{align*}
\sigma: &= 2 \cdot {p-3 \over 2} - \sum(a+ \bar a) = p-3 - \sum^{p-2}_{i=2} i 
\\
&= p-3+1 - 
{(p-2)(p-1) \over 2} \equiv -2-{(-2)(-1) \over 2} 
=- 3 \ (\bmod \ p). 
\end{align*}
Thus,
\begin{equation*}
\Pi^{p-2}_{a=2} [a]_{{\dps{q}}^{\bar a-1}} \equiv q^{-3} (\bmod \ [p]). \tag{132a}
\end{equation*}
In addition, 
\begin{equation*}
[1]_{{\dps{q}}^{0}} = 1 = q^0, \tag{132b}
\end{equation*}
and
\begin{equation*}
[p-1]_{{\dps{q}}^{p-1-1}} \equiv [p-1]_{{\dps{q}}^{-2}} = {1-q^{-2(p-1)} \over 1-q^{-2}} \equiv 
{1 - q^2 \over 1-q^{-2}} = -q^2 \ (\bmod \ [p]). \tag{132c}
\end{equation*}
Combining formulae (132a)--(132c) we obtain (126a).  Since 
\begin{equation*}
[p-1] = [p] - q^{p-1} \equiv -q^{-1} \ (\bmod \ [p]), \tag{133}
\end{equation*}
(126b) follows from (126a).
\end{proof}
\begin{remark}{134}
Suppose $p \equiv 1$ (mod 4) and $a$ is such that $0 < a < p$ and 
\begin{equation*}
a^2 \equiv -1 \ (\bmod \ p). \tag{135a}
\end{equation*}
Then $\bar a = p - a$, because
\begin{equation*}
a(p - a) \equiv 1 \ (\bmod \ p), \ \ \ 0 < \bar a < p. 
\end{equation*}
Therefore, by formula (131), 
\begin{equation*}
[a]_{q^{-a-1}} [-a]_{q^{a-1}} \equiv q^2 \ (\bmod \ [p]). \tag{135b}
\end{equation*}
Since
\begin{equation*}
[-a]_Q = - Q^{-a} [a]_Q, 
\end{equation*}
formula (135b) can be rewritten as 
\begin{equation*}
[a]_{q^{-1-a}} [a]_{q^{-1+a}} \equiv - q^{a(a-1)+2} \equiv q^{(1-a)} \ (\bmod \ [p]), \tag{135c}
\end{equation*}
a quantum version of formula (135a). 
\end{remark}

We see that QNT subsumes the classical cyclotomy.  But we have barely scratched the surface 
of QNT.  There remain plenty of eminently mentionable 
but unmentioned in this review classical formulae, theorems, and arguments where powers of 
integers enter $-$ and many such can be found in {\it{Ramanujan's Notebooks}}; whether all 
these classical results can be quantized or not,  no one can know in advance $-$ the gods 
do not respond to mathematical queries, and Ramanujan is no longer with us.  But the 
quantization country is open to all willing to explore it, and while success is not guaranteed, 
adventure, excitement, and bewilderment are.

To illustrate the latter, let's conclude this review with the following example.  In classical 
theory, the SFT:  ``$a^{p-1}\equiv 1 \ (\bmod \ p)$ for {\it{a}} coprime to $p$'' can be used 
directly to deduce the Wilson theorem ``$(p-1)!\equiv -1 \ (\bmod \ p)$'', as follows.  The 
congruence
\begin{equation*}
x^{p-1} - 1 \equiv 0 \ (\bmod \ p) \tag{136}
\end{equation*}
has $p-1$ solutions
\begin{equation*}
x = 1, 2, ... , p-1. \tag{137}
\end{equation*}
Therefore, 
\begin{equation*}
x^{p-1} - 1 \equiv (x-1) (x-2) ... (x-(p-1)) \ (\bmod \ p) \tag{138}
\end{equation*}
as polynomials in $x$, and the $x$-free terms in these polynomials yield:
\begin{equation*}
-1 \equiv (-1) (-2) ... (-(p-1)) = (p - 1)! \ (\bmod \ p). \tag{139}
\end{equation*}
Let us apply this argument to our $q$-version (96) of SFT:  
\begin{gather*}
< a^{p-1}> \equiv 1 \ (\bmod \ [p]), \ \ \ \forall a = 1, ... , p - 1 \ \ \Leftrightarrow 
\tag{140a}
\\[1ex]
{\Pi^{p-1}_{k=1} (1-q^{kx}) \over \Pi^{p-1}_{k=1} (1 - q^k)} \equiv 1 \ 
(\bmod \ [p]), \ \ \ \forall x=1, 2, ..., p -1 \ \Leftrightarrow \tag{140b}
\\[1ex]
\Pi^{p-1}_{k=1} (1 - q^{kx}) \equiv \Pi^{p-1}_{k=1} (1-q^k) \ {\rm{[by}} \ 
(108)] \equiv p \ (\bmod \ [p]) \ \Leftrightarrow \tag{140c}
\\[1ex]
\Pi^{p-1}_{k=1} (y^k -1)  -p \equiv 0 \ (\bmod \ [p]), \ \ \ \forall y \ (=q^x) = q^1, ... , 
q^{p-1} \Leftrightarrow \tag{140d}
\\[1ex]
 \Pi^{p-1}_{k=1} (y^k -1) -p \equiv \Pi^{p-1}_{k=1} (y-q^k) \ (\bmod \ [p]). 
\tag{141}
\end{gather*}
Now, 
\begin{equation*}
\Pi^{p-1}_{k=1} (y - q^k) \equiv 1+y + ... + y^{p-1} = {y^p -1 \over y-1} \ (\bmod \ [p]), 
\tag{142}
\end{equation*}
because both sides vanish $\bmod \ (q^p-1)$ whenever $y$ is a primitive $p$-th root of unity.  
Thus, the congruence (141) can be rewritten as
\begin{equation*}
\Pi^{p-1}_{k=1} (y^k-1) - p \equiv 1+y + ... + y^{p-1} = {y^p -1 \over y-1} \ (\bmod \ [p]). 
\tag{143}
\end{equation*}
But there are hardly any $q$'s left in this congruence, so we conclude that the polynomial 
\begin{equation*}
\varphi_p (y) = \Pi^{p-1}_{k=1} (y^k -1) -p \tag{144a}
\end{equation*}
is divisible by the polynomial 
\begin{equation*}
\psi_p (y) = 1 + y + ... + y^{p-1} = {y^p-1 \over y-1} . \tag{144b}
\end{equation*} 
The ratio
\begin{equation*}
\chi_p (y) = {\varphi_p (y) \over \psi_p (y)} = (\Pi^{p-1}_{k=1} (y^k -1) -p) (y-1)/(y^p-1) 
\tag{145}
\end{equation*}
is a polynomial in $y$ of degree $d_p = {({p \atop 2})} - (p-1) = {{({p -1 \atop 2})}}$ .  For 
$p=3, \ d_3 = 1$, and 
\begin{equation*}
\chi_3 (y) = y-2 \ \ \ \ ( \Rightarrow x_{root} = \log_q 2). \tag{146}
\end{equation*}
For $p=5, \ d_5 = 6$, and
\begin{equation*}
\chi_5 (y) = y^6 - 2y^5 + y^3 + 3y -4. \tag{147}
\end{equation*}
What are the roots of these polynomials, and what is their meaning? ... Silence.  Another 
mystery...

This  part of the review of Brendt's magnificient masterpiece can be  summarized as 
\begin{theorem}{148}
Reading \emph{Ramanujan's Notebooks} could be 
 good for one's mathematical health.
\end{theorem}

\newpage
\section*{References}
\begin{flushright}
\parbox{5.5cm}{\sffamily \small `Always verify references'.  
This is so absurd in mathematics that I used to say 
provocatively:  `never...' \\
\null\qquad J. E. Littlewood}
\end{flushright}

\renewcommand{\refname}{}

\strut\hfill

\noindent
$\bar{\mbox{S}}$.D.G.
\label{lastpage}

\vskip 12pt
\centerline{Boris A Kupershmidt}

\end{document}